\documentclass[12pt]{amsart}
\usepackage[active]{srcltx}

\usepackage{amsfonts,amssymb, enumerate}

\usepackage[utf8]{inputenc}  
\usepackage[T1]{fontenc}

\def\Pr{\begin{proof}}
\def\Rp{\end{proof}}

\def\IR{\mathbb R}

\def\IN{\mathbb N}

\def\IQ{\mathbb Q}

\def\IC{\mathbb C}
\def\IF{\mathbb K}
\def\IF{\mathbb F}
\def\IZ{\mathbb Z}
\def\IL{\mathbb L}
\def\IF{\mathbb K}
\def\IF{\mathbb F}

\def\ZF{\mathbf {ZF}}
\def\ZFA{\mathbf {ZFA}}

\def\DC{\mathbf {DC}}
\def\HB{\mathbf {HB}}
\def\BPI{\mathbf {BPI}}
\def\AC{\mathbf {AC}}
\def\VKM{\mathbf {VKM}}
\def\KM{\mathbf {KM}}
\def\MC{\mathbf {MC}}
\def\RL{\mathbf {RL}}

\def\I{\mathbf {I}}

\DeclareMathOperator{\spanv}{span}

\DeclareMathOperator{\diam}{diam}

\DeclareMathOperator{\Prime}{PRIME}
\DeclareMathOperator{\rk}{rk}

\DeclareMathOperator{\On}{On}

\theoremstyle{plain}

\newtheorem{corollary}{Corollary}
\newtheorem{proposition}{Proposition}
\newtheorem*{proposition*}{Proposition}
\newtheorem{theorem}{Theorem}
\newtheorem*{theorem*}{Theorem}
\newtheorem*{corollary*}{Corollary}
\newtheorem{lemma}{Lemma}
\newtheorem*{lemma*}{Lemma}

\theoremstyle{definition}

\theoremstyle{remark}
\newtheorem{remark}{Remark}

\newtheorem{example}{Example}

\date{\today}
\usepackage{graphicx}

\begin{document}
\title{Multiple Choices imply the Ingleton and Krein-Milman axioms}
\author[M.~Morillon]{Marianne Morillon}
 \address{Laboratoire d'Informatique et Mathématiques, 
Parc Technologique Universitaire, Bâtiment 2,
2 rue Joseph Wetzell, 97490 Sainte Clotilde}
 \email{Marianne.Morillon@univ-reunion.fr}

 \begin{abstract} In set theory without the Axiom of Choice, we consider
 Ingleton's axiom which is  the counterpart in ultrametric analysis of the Hahn-Banach axiom. We show that in $\ZFA$, set theory without the Axiom of Choice weakened to allow ``atoms'',  Ingleton's axiom does not imply the Axiom of Choice (this solves in $\ZFA$ a question raised by van Rooij, \cite{Rooij}). 
  We also prove that in $\ZFA$, the ``multiple Choice'' axiom implies the Krein-Milman axiom. We deduce that, in $\ZFA$, the conjunction of the Hahn-Banach, Ingleton and Krein-Milman axioms  does not imply the Axiom of Choice. 
 \end{abstract}

 \maketitle

\section{Introduction}
We denote by $\ZFA$ (see \cite[p.~44]{Je73}) the set theory $\ZF$ with the axiom of  extensionality  weakened to allow the existence of atoms. 
We denote by $\MC$ (``Multiple Choice'')  the following consequence of the Axiom of Choice $\AC$ (see \cite[p.~133]{Je73} and
 form~67 of \cite[p.~35]{Ho-Ru}): 
{\em ``For every infinite family $(X_i)_{i \in I}$ of nonempty sets, there exists a family 
$(F_i)_{i \in I}$ of nonempty finite sets such that for each $i \in I$, $F_i \subseteq X_i$.''}

It is known that $\MC$ is equivalent to $\AC$ in $\ZF$ (see \cite[Theorem~9.1 p.~133]{Je73}), but $\MC$ does not imply $\AC$ in $\ZFA$ (see \cite[Theorem~9.2 p.~134]{Je73}).
Given a prime number $p \ge 2$, we also consider the following refined  statement for each  prime number  $p \ge 2$, where for every {\em finite} set $F$, we denote by $\#F$ the cardinal of $F$:\\
 $\MC(p)$: {\em ``For every family $(A_i)_{i \in I}$ of nonempty sets, there exists a family $(B_i)_{i \in I}$ of finite  sets such that for every  $i \in I$, $B_i \subseteq A_i$ and  $\#B_i$ is not a multiple of  $p$.''} 
 
 The conjunction of all statements $\MC(p)$ for all prime numbers $p$ is denoted by form~218 in \cite[p.~52]{Ho-Ru}).
 Levy (1962) built a model $\mathcal N6$ of $\ZFA$ (see \cite[p.~185]{Ho-Ru}, \cite{Levy62}) satisfying $\MC(p)$ for every prime number $p \ge 2$ (and thus satisfying also $\MC$), in which there exists a sequence   $(F_n)_{n \in \IN}$ of finite sets such that for every $n \in \IN$,  $\#F_n = n+1$ and  $\prod_{n \in \IN}F_n= \varnothing$: such a model  does not satisfy  $\AC$ (and more precisely, this model does not satisfy the ``countable Axiom of Choice for finite sets''). It follows that in 
 $\ZFA$, the conjunction of  statements $\MC(p)$ for all prime numbers $p$  does not imply $\AC$.

\medskip

Consider   the well-known Hahn-Banach statement $\HB$ (form~52 of \cite{Ho-Ru}):\\
$\HB$: {\em Given a real vector space $E$ and a sublinear functional $p:E \to \IR$, there exists a linear functional $f:E \to \IR$ such that $f \le p$.}

It is known  (see \cite[Corollary~1]{Mo-Rado}) that in $\ZFA$,  $\MC$ implies Rado's selection Lemma $\RL$ (form~99 of \cite{Ho-Ru}),
and that  $\RL$  implies $\HB$ (see \cite[Theorem~1]{Mo-Rado}), thus $\MC$ implies $\HB$.

 Given a spherically complete ultrametric  valued field $(\IF,|.|)$ 
 (see Section~\ref{subsec:sph-comp}), Ingleton proved 
 in set theory $(\ZF \mathbf + \AC)$  (see \cite{Ing52}), for each ultrametric normed space over a ``spherically complete'' ultrametric valued field $(\IF,|.|)$, a ``Hahn-Banach''-type result which we denote by 
 $\I_{(\IF,|.|)}$. 
 A.C.M. van Rooij  (\cite{Rooij}) showed that for each ultrametric valued field 
 $(\IF,|.|)$ such that the large ball $\{x \in \IF : |x| \le 1\}$ of $\IF$ is compact (whence $(\IF,|.|)$ is spherically complete), $\BPI$ implies $\I_{(\IF,|.|)}$. He asked whether  the  ``full Ingleton theorem'' ({\em i.e.} the conjunction of all statements $\I_{\IF,|.|}$ for all  spherically complete  ultrametric valued fields  $(\IF,|.|)$) implies $\AC$. We shall  show that in set theory $\ZFA$ (set theory without choice weakened to allow ``atoms''), the ``full Ingleton theorem''  does not imply $\AC$. 
 More precisely, for every prime number $p \ge 2$, 
$\MC(p)$ implies the 
statement $\I_{(\IF,|.|)}$ for each  spherically complete ultrametric valued field $(\IF,|.|)$ with null characteristic such that the restriction of the absolute value on $\IQ$ is  equivalent to the $p$-adic absolute value, and  also for  each spherically complete ultrametric valued field $(\IF,|.|)$ with  characteristic $p$   (see Corollary~\ref{coro:MCptoIng}-\eqref{it:MCptoIngp} in Section~\ref{sec:mcp_to_ing}).
In $\ZFA$, $\MC$ implies the 
statement $\I_{(\IF,|.|)}$ for each spherically complete ultrametric valued field $(\IF,|.|)$ with null characteristic such that the restriction of the absolute value on $\IQ$ is trivial (see Corollary~\ref{coro:MCptoIng}-\eqref{it:MCtoIng0} in Section~\ref{sec:mcp_to_ing}). 
Since Levy's above model of $\ZFA$ satisfies  $\MC(p)$ for every prime number $p \ge 2$ (and thus it also satisfies $\MC$), it follows that 
in $\ZFA$, the ``full Ingleton axiom'' does not imply $\AC$. 

In  Section~\ref{sec:MC2KM}, we prove (see Theorem~\ref{theo:mc2hb}) that in $\ZFA$, 
$\MC$ implies the Krein-Milman statement  (form~65 of \cite{Ho-Ru}):
\par $\KM$: {\em Given a Hausdorff locally convex topological real vector space $E$, every nonempty compact convex subset of $E$ has an extreme point.}\\
We shall use  ``trees'' of ``facets'' of a  nonempty convex compact subset of a topological real vector space. We deduce that 
$\MC \Rightarrow (\HB \mathbf + \KM)$. This generalizes  
Pincus's following result (\cite{Pi72b}):  the model $\mathcal N_2$  (described in  \cite[p.~178]{Ho-Ru}) of  
$\ZFA \mathbf + \neg \AC \mathbf + \MC$  satisfies $(\KM \mathbf + \HB)$. 
It follows that the conjunction of statements $\MC(p)$ for every prime number $p \ge 2$ implies Rado's lemma (and hence $\HB$), the ``full Ingleton axiom'' and also the Krein-Milman axiom. 
In particular,   Levy's model $\mathcal N6$ of $\ZFA \mathbf + \neg \AC$ satisfies all these statements.

The paper is organized as follows: in Section~\ref{subsec:def-not}, we sketch basic results about ultrametric valued fields, in Section~\ref{sec:mcp_to_ing} we prove in $\ZFA$ that the conjunction  of the statements $\MC(p)$ for all prime numbers $p$ implies the ``full Ingleton Theorem'', and in the last Section~\ref{sec:MC2KM}, we prove in $\ZFA$ that $\MC$ implies the ``Krein-Milman'' statement. 

\section{Preliminaries on ultrametric valued fields} \label{subsec:def-not}
\subsection{Valued fields}
An {\em absolute value} on a (commutative) field $\IF$ is a mapping $|.|: \IF \to \IR_+$  satisfying the following three properties for every $\lambda$, $\mu \in \IF$: $|\lambda|=0 \Leftrightarrow \lambda=0$; $|\lambda \mu|=|\lambda| |\mu|$ and
 $|\lambda + \mu| \le |\lambda|+|\mu|$.  Each valued field $(\IF,|.|)$ gives rise to a metric $d_{|.|}: \IF \times \IF \to \IR_+$ defined by $d(x,y)=|x-y|$ for every $x,y \in \IF$.  An absolute value $|.|$ on $\IF$ is said to be {\em ultrametric} if the associated metric $d_{|.|}$ is ultrametric ({\em i.e.} $\forall x,y,z \in \IF \; d(x,z) \le \max(d(x,y),d(y,z))$), equivalently if  
 for every $\lambda$, $\mu \in \IF$, $|\lambda + \mu| \le \max(|\lambda|, |\mu|)$.

 \begin{example}[The trivial absolute value on a field]
 For each commutative field $\IF$, the mapping $|.|_{triv} : \IF \to \IR_+$ associating to each $\lambda \in \IF$ the real
   number $0$ if $\lambda=0$ and $1$ otherwise is an absolute value, called the {\em trivial} absolute value on $\IF$. The associated metric is the discrete distance $d_{disc}$ on $\IF$, satisfying  
   $d_{disc}(x,y)=0$ if $x = y$ and $1$ otherwise. 
\end{example}

 \begin{remark} Notice that if $|.|$ is  an absolute value on a field $\IF$, then $|1_{\IF}|=1$. 
 If $\IF$ is a finite field with $m$ elements, then for every $x \in \IF \backslash \{0\}$, 
 $x^{m-1}=1$ thus $|x|^{m-1}=1$ so $|x|=1$; it follows that the trivial absolute value is the only  absolute  value  on $\IF$.
\end{remark}
 
 \subsection{Ultrametric valued fields}
 Given a non-empty set $X$ and a semi-metric $d: X \times X \to \IR_+$,  a {\em large} ball of the semi-metric space $(X,d)$ is a ball with  ``large inequalities''  {\em i.e.} a subset of $X$  of the form  
$B(a,r) := \{x \in X : d(x,a)  \le r\}$ where $a \in X$ and $r \in \IR_ +$.
  Given a ultrametric valued field $(\IF,|.|)$, the large ball $B_{\IF}(0,1)$ is a subring of $\IF$, and the ``strict ball'' 
  $\mathfrak M_{\IF} :=\{x \in \IF : |x|<1\}$ is the unique maximal ideal of the ring $B_{\IF}(0,1)$ ($\mathfrak M_{\IF}$ is the set of non-invertible elements of the ring $B_{\IF}(0,1)$). The field $B_{\IF}(0,1)/\mathfrak M_{\IF} $ is the {\em residue class field} of the valued field $\IF$. 
  The mapping $|.|_{\restriction \IF \backslash \{0\}} : \IF \backslash \{0\} \to {\IR^*_+}$ is a group morphism {\em w.r.t.} the multiplicative laws, thus the set $V_{(\IF,|.|)}:=\{|x| : x \in \IF \backslash \{0\}\}$ is a subgroup 
  of the multiplicative group ${\IR^*_+}$, which is called the {\em value group} of the absolute value $|.|$.
  The value group $V_{(\IF,|.|)}$ of a ultrametric absolute value $|.|$ on a valued field $\IF$ is either discrete or dense in $\IR^ *_+$.
  
  \begin{example}[The trivial absolute value on a field]
 For each commutative field $\IF$ endowed with the trivial absolute value $|.|_{triv} : \IF \to \IR_+$,
then   $B_{\IF}(0,1)=\IF$, and the maximal ideal $\mathfrak M_{\IF}$ is the null ideal $\{0\}$, thus the residue class field of $\IF$ is $\IF$. The group of values $V_{\IF}$ is the trivial subgroup $\{1\}$ of
$\IR^*_+$.
\end{example}

\begin{example}
Given a ultrametric valued field $(\IF,|.|)$ such that 
 $V_{(\IF,|.|)}$  is discrete and $V_{\IF} \neq \{1\}$, then  $s:= \max \{|x| : x \in \mathfrak M_{\IF}\} \in ]0,1[$, and   $V_{(\IF,|.|)}= \{s^n : n \in \IZ\}$. 
\end{example}

\subsection{Cauchy-complete semi-metric spaces}
Given a semi-metric space $(X,d)$, a filter $\mathcal F$ on $X$ is {\em Cauchy} 
if for every real number $\varepsilon >0$, there exists $F \in \mathcal F$ such that the {\em $d$-diameter} 
$\sup_{x,y \in F} d(x,y)$ of $F$ is $< \varepsilon$.  
 A semi-metric space $(X,d)$ is {\em Cauchy complete} if every Cauchy filter of this semi-metric space 
has a non-empty intersection. 

\begin{example}
\begin{enumerate}
  \item The field $\IR$ of real numbers endowed with the usual absolute value is Cauchy complete.
For every set $X$, denoting by $B(X,\IR)$ the vector space of bounded functions $f: X \to \IR$, 
then the normed vector space $B(X,\IR)$ endowed with the uniform norm is Cauchy-complete.
  \item Every set $X$  endowed with the discrete distance $d_{disc}$ is Cauchy-complete since, for every real number 
$\varepsilon \in ]0,1[$, non-empty $\varepsilon$-small subsets of $X$ are singletons. In particular, 
every valued field endowed with the trivial absolute value is Cauchy-complete. 
\end{enumerate}
\end{example}

\begin{remark} Say that a semi-metric space $(X,d)$  is {\em sequentially complete} if every $d$-Cauchy 
sequence of $X$ converges to a point of $X$.  There are models of $\ZF$ (for example the ``basic Cohen model'' described in \cite[Section~5.3]{Je73}) in which there exists a dense (hence infinite) subset $D$ of $\IR$ which is {\em Dedekind-finite} ({\em i.e.} which does not contain any subset equipotent with $\IN$). Such a set $D$ endowed with the usual distance of $\IR$ is a metric space which is ``sequentially complete'' (because for each sequence $(x_n)_{n \in \IN}$ of $D$, the set $\{x_n: n \in \IN\}$ is finite thus each Cauchy-sequence of $(D,d)$ is stationnary) but  the metric subspace $(D,d)$ of $\IR$ is not Cauchy-complete since $D$ is not closed in $\IR$. 
\end{remark}

A {\em Cauchy-completion} of a semi-metric space $(X,d)$ is a Cauchy complete metric space  $(\hat X,\hat d)$ 
together with an isometry $j: (X,d) \to (\hat X,\hat d)$  such that $j[X]$ is dense in $\hat X$.

\begin{proposition} For every semi-metric space $(X,d)$, the mapping $j: X \to B(X,\IR)$ 
associating to each $a \in X$ the mapping $f_a=(d(x,a))_{x \in X}$ is isometric, the closed subspace
$\overline{j[X]}$ of $B(X,\IR)$  is Cauchy-complete, and  $(j,\overline{j[X]})$ is a Cauchy-completion 
of $(X,d)$.
\end{proposition}

Given a valued field $(\IF,|.|)$, the laws $+$ and $.$ of the field $\IF$ are uniquely extendable into continuous laws on the Cauchy-completion $\hat \IF$ of the metric space $(\IF,d_{|.|})$ (the law $+$ is uniformly continuous and the law $.$ is uniformly continuous on  $B_{\IF}(0,1) \times B_{\IF}(0,1)$). The structure  $\hat \IF$ 
endowed with these extended laws is a valued field.

 \subsection{Non trivial absolute values on $\IQ$}
 \subsubsection{Archimedean absolute values on $\IQ$}
 An absolute value on a field $\IF$ is {\em archimedean} if  
 $\sup \{|n.1_{\IF}| : \;  n \in \IN\}=+\infty$. Given a field $\IF$, 
 an absolute  value on $\IF$ is non archimedean if and only if it is ultrametric (see \cite[Theorem~1.1 p.~2]{Rooij78}). 
  
\begin{example}
 \begin{enumerate}
 \item The mapping $|.|_{\IQ} : \IQ \to \IR_+$ associating to each rational number $x$ le number $\max(x,-x)$
 is an archimedean absolute value on the field $\IQ$ of rational numbers.
 \item Denoting by $\IR$ the ordered field of real numbers, {\em i.e.} the Cauchy-completion of the valued field  $(\IQ, |.|_{\IQ})$, the   ``usual'' absolute value   $|.|_{\IR}: \IR \to \IR_+$ which associates to each $x \in \IR$  the real number $\max(x,-x)$ is archimedean.  
 \item  Denoting by  $\IC$ the field of complex numbers, the ``modulus function'' 
 $|.|_{\IC} : \IC \to \IR_+$ is an archimedean absolute value on $\IC$. 
 \item  Given any subfield  $\IF$ of $\IC$, the restriction of the ``modulus function'' 
 ${|.|_{\IF}}:={|.|_{\IC}}_{\restriction \IF} : \IF \to \IR_+$ is an archimedean absolute value on $\IF$. 
 \end{enumerate}
 \end{example}
  
 \begin{remark} Notice that,   given a subfield $\IF$ of $\IC$, for each real number $\tau \in ]0,1]$, 
 the function ${{|.|_{\IC}}^{\tau}}_{\restriction \IF}$ is also an archimedean absolute value on $\IF$. 
  Given a ultrametric field $(\IF, |.|)$, for each real number $\tau >0$, 
 the function ${|.|}^{\tau}$ is also a ultrametric absolute value on $\IF$  which is equivalent to the
 initial absolute value $|.|$.  
 \end{remark}

\subsubsection{$p$-adic absolute values on $\IQ$}
We denote by $\Prime$ the set of prime natural numbers. 
For each $p \in \Prime$, we denote by   $v_p: \IQ \to \IZ \cup \{+\infty\}$ 
the  {\em $p$-adic valuation} on $\IQ$, such that $v_p(0)=+\infty$, $v_p(1)=v_p(-1)=0$, and for all $x \in \IQ \backslash \{0,1,-1\}$, $v_p(x)$ is the exponent $\alpha$
of $p$ in the prime decomposition $x=\pm p^{\alpha} \prod_{q \in \Prime; q \neq p}q^{\alpha_q}$ 
of $x$ in prime numbers. 
For each prime natural number $p$, we denote by $\IF_p=\{0,1,\dots,p-1\}$ the finite field with $p$ elements. 

\begin{proposition}[{\cite[Prop.~2.4.3 p.~40]{Gouv}}]
\begin{enumerate}
\item \label{it:val_Qp1} The mapping $x \mapsto |x|_p:=p^{-v_p(x)}$ is a ultrametric  absolute value on the field $\IQ$ of rational numbers.
\item \label{it:val_Qp2} The  group of values of  of the valued field $(\IQ,|.|_p)$ is $\{p^n : n \in \IZ\}$; 
\item \label{it:subring_Qp} The subring $B(0,1)$ of $(\IQ,|.|_p)$ is the subset 
$\{\frac{m}{n} : m \in \IZ; n \in \IN^* \text{ s.t. } p \nmid n \}$ and thus $\IZ$ is a subring of $B(0,1)$.
\item \label{it:max_ideal} The maximal ideal $\mathfrak M_{(\IQ,|.|_p)}$ is $pB(0,1)$.
\item \label{it:cosets_residue_class} The {\em cosets} of the ring $B(0,1)$ {\em modulo} the maximal ideal  $pB(0,1)$ are the sets $i+pB(0,1)$ where $i \in \{0, \dots, p-1\}$, and  the residue class field is the finite field $\IF_p$.  
\end{enumerate}
\end{proposition}
\Pr \eqref{it:val_Qp1} and \eqref{it:val_Qp2} are trivial. See the proofs of \eqref{it:subring_Qp} and  \eqref{it:max_ideal} in 
 \cite[Prop.~2.4.3-(i) and (ii) p.~40]{Gouv}. 
\eqref{it:cosets_residue_class} We first notice that if $i,j \in \{0,\dots,p-1\}$, then the sum of the two classes 
$i+pB(0,1)$ and $j+pB(0,1)$ is the class $r+pB(0,1)$ where $r$ is 
the remainder in the euclidean division of $(i+j)$ by $p$: indeed, let $q \in \IN$ 
such that $i+j=qp+r$ where $r \in \{0,\dots,p-1\}$; 
then $|i+j-r|_p = |qp|_p$ which is $0$ if $q=0$ and $p^{-\alpha}$ where $\alpha \in \IN^*$ else, thus 
$|i+j-r|_p <1$ so $(i+j)+pB(0,1)=r+pB(0,1)$. 
Let $x \in \IQ$ such that $|x|_p=1$. Let us show that there exists $i \in \{1, \dots, p-1\}$ such that 
$|x-i|_p<1$. It is sufficient to prove this if $x= \frac{1}{n}$ where  $n \in \IN^*$ with $|n|_p=1$ 
{\em i.e.} $p$ does not divide $n$. Since $p \nmid n$, $n+p\IZ$ is invertible in $\IZ/p\IZ$, thus 
let $i \in \{1,\dots,p-1\}$ such that 
$ p/ (ni -1)$; then $|n i-1|_p <1$ {\em i.e.} $|n|_p |i - \frac{1}{n}|_p < 1$ {\em i.e.} 
$|i - \frac{1}{n}|_p < 1$ thus $\frac{1}{n} -i \in pB(0,1)$. 
\Rp

\begin{proposition}[{\cite[Lemma~3.2.3 p.~51]{Gouv}}]
\begin{enumerate}
\item \label{it:IQ_not_complete} The valued field $(\IQ,|.|_p)$ is not complete.
\item Denoting by $\IQ_p$ the  Cauchy-completion of the valued field  $(\IQ, |.|_p)$, the  extended  absolute value   $|.|_{p}: \IQ_p \to \IR_+$ is also ultrametric and has the same group of values as $(\IQ, |.|_p)$.  Moreover, $pB_{\IQ_p}(0,1)$ is still a maximal ideal of  the subring $B_{\IQ_p}(0,1)$ of $\IQ_p$
and the quotient field  $B_{\IQ_p}(0,1)/pB_{\IQ_p}(0,1)$ is still the finite field $\IF_p$.
\end{enumerate}
\end{proposition}

\begin{theorem*}[Ostrowski's theorem, {\cite[Th.~ 1.2 p.~3]{Rooij78}}]
 Every  non trivial absolute value on $\IQ$ is of the form $|.|_{\IR}^{\tau}$ where $\tau$ is some real number such that $0<\tau<1$, or of the form $|.|_p^{\tau}$  for some  prime number $p$ and  $\tau >0$.
 \end{theorem*}

\subsection{Spherically complete ultrametric spaces} \label{subsec:sph-comp}
A ultrametric semi-metric space  is {\em spherically complete} if every chain 
of large balls  of this semi-metric  space has a non-empty intersection.

\begin{example} 
  Each nonempty set $X$ endowed with the discrete metric $d_{disc}$  is spherically complete, because large balls of the metric space $(X,d_{disc})$  are singletons of $X$ or the whole set  $X$. 
  \end{example}

\begin{proposition}
Every spherically complete ultrametric space is Cauchy-complete.
\end{proposition}
\Pr Assume that  $(X,d)$ is a ultrametric metric space which is spherically complete and that 
$\mathcal F$ is a Cauchy filter of the metric space $(X,d)$. 
Consider the set $\mathcal C$ of large balls which belong to $\mathcal F$. 
For every $B_1$, $B_2 \in \mathcal C$, $B_1 \cap B_2 \in \mathcal F$ thus $B_1 \cap B_2 \neq \varnothing$;
since $d$ is ultrametric, it follows that
 $B_1 \subseteq B_2$ or $B_2 \subseteq B_1$. 
 Thus $\mathcal C$ is a chain of large balls. Notice that for every $F \in \mathcal F$ with diameter 
 $d_F$, for every $x \in F$, $F \subseteq B(x,d_F)$, thus, 
 for every $x, y \in F$, $B(x,d_F)=B(y,d_F)$.  Since $(X,d)$ is spherically complete, 
 the set $\cap \mathcal C$ is nonempty. 
Since $\mathcal F$ is a Cauchy-filter, $\mathcal C$ contains  balls of  diameter $\varepsilon$ for every 
$\varepsilon >0$, thus  $\cap \mathcal C$ is a singleton
$\{a\}$. We now show that for every $\varepsilon >0$ and every $F \in \mathcal F$, 
$B(a,\varepsilon) \cap F \neq \varnothing$ (whence it will follow that $a \in \overline{F}$). It is sufficient to prove that for every $\varepsilon >0$ and every $F \in \mathcal F$ such that $\diam(F)< \varepsilon/2$, $B(a,\varepsilon) \cap F \neq \varnothing$.
Given some $F \in \mathcal F$ such that $\diam(F)< \varepsilon /2$, then, for every $x \in F$, $F \subseteq B(x,\varepsilon/2)$, thus 
$B(x,\varepsilon/2) \in \mathcal C$, so $a \in B(x, \varepsilon/2)$ whence $B(a, \varepsilon/2) \cap F \neq \varnothing$. 
\Rp

\begin{remark}
For each prime number $p$, the valued field $(\IQ,|.|_p)$ is not complete hence not spherically complete,
  \end{remark}

\subsection{For each prime number $p$, $\IQ_p$ is locally compact hence spherically complete}  
\subsubsection{The topological ring $\IZ_p$}
Let $p$ be a  prime number. For  each natural number $n$, we endow the finite ring $\IZ/{p^n\IZ}$ with the discrete topology 
and we endow the product ring   $\mathcal A:=\prod_{n \in \IN^*} \IZ/p^n\IZ$ with the product topology: thus  
the topological ring  $\mathcal A$ is compact and  metrizable (by the usual metrics on the product of metric spaces 
$\mathcal A$). We denote by $\IZ_p$  (see \cite[Prop.~1.2.1 p.~18]{Amice75}) the closed subring  
$\{x=(x_n)_{n \in \IN^*} \in \mathcal A : \forall n \in \IN^* \;\;  p^n \mid (x_{n+1}-x_n)\}$ of $\mathcal A$.
It follows that $(\IZ_p,+,\times)$ is a commutative compact topological ring. 
The topological subspace $\IZ_p$ of $\mathcal A$ is also metrizable.
 The unit of $\IZ_p$ is  the constant sequence $(1)$. 
Let $can: \IZ \to \IZ_p$ be the canonical injection. 
 If $a \in \IN$, if $a=\sum_{k}a_kp^k$ is the $p$-ary expansion of $a$ (where each $a_k \in \{0,\dots,p-1\}$,
then $can(a)$ is the element $x=(x_n)_{n \in \IN^*} $ of $\IZ_p$ such that for each $n \ge 1$, $x_n=\sum_{0 \le k <n}a_kp^k$.
It follows that  $can : \IZ \to \IZ_p$ is dense  (see \cite[Prop.~1.2.3 p.~19]{Amice75}.

\begin{proposition}[{\cite[Prop.~1.4.5. p.~23]{Amice75}}]
An element $x=(x_n)_{n \in \IN^*}$ is invertible in $\IZ_p$ if and only if $x_1 \neq 0$ in $\IZ/p\IZ$.
\end{proposition}

 \subsubsection {The topological ring $\IZ_p$ and the large ball $B(0,1)$ of $\IQ_p$ are isomorphic}
 \begin{proposition}
The topological ring  $\IZ_p$ and the large unit ball  $B(0,1)$ of $\IQ_p$ are isomorphic.
It follows that the large ball $B(0,1)$ of $\IQ_p$ is compact. 
\end{proposition} 
\Pr For each $x=(x_n)_{n \in \IN^*}  \in  \IZ_p$, considering each $x_n \in \IZ/{p^n\IZ}$ as an element of $\{0,\dots,p^n-1\} \subseteq \IZ$,
 the sequence  $(x_n)_{n \ge 1}$ of $\IQ$ is a Cauchy sequence for the metric $d_p$ on $\IQ$, thus it converges to a unique element  of $B_{\IQ_p}(0,1)$.
 Let $f: \IZ_p \to B_{\IQ_p}(0,1)$ be the mapping associating to each $x=(x_n)_{n \in \IN^*}  \in  \IZ_p$ its limit in $B_{\IQ_p}(0,1)$.
Then $f$ is a morphism of rings, and $f$ is one-to-one. 
Moreover, the  mapping $f$ is onto: since  $f[\IZ_p]$ is  a compact (hence closed) 
subset of $B(0,1)$, it is sufficient to prove that 
$(\IQ \cap B_{\IQ_p}(0,1)) \subseteq f[\IZ_p]$, or, equivalently, that for  each $n \in \IN^*$ such that  $p \nmid  n$, 
$\frac{1}{n} \in f[\IZ_p]$: given $n \in \IN^*$ such that  $p \nmid  n$, consider the $p$-ary expansion $n=\sum_{k}a_kp^k$
where each $a_k \in \{0,\dots,p-1\}$, then $a_0 \neq 0$ thus $can(n)=(a_0, a_0+pa_1, \dots)$ is invertible in $\IZ_p$, so $\frac{1}{n} \in f[\IZ_p]$.
\Rp

  \begin{corollary}
The Cauchy completion $\IQ_p$ of $(\IQ,|.|_p)$ is  spherically complete.
  \end{corollary}
  \Pr Using the previous Proposition, every large ball  of $\IQ_p$ is compact whence the metric space 
  $(\IQ,|.|_p)$ is Cauchy-complete.
  \Rp

 \begin{remark} Since the field  $\IQ$ is countable, $\IQ$ has a unique countable algebraic closure $\IQ_{ac}$. 
 For each prime number $p$, there exists a unique absolute value on $\IQ_{ac}$ extending $|.|_p$ 
 (using \cite[Prop.~2.6.1 p.~63]{Amice75}). 
 The Cauchy-completion of the valued field $(\IQ_{ac},|.|_p)$ is still algebraically closed and is denoted by $\IC_p$: thus 
  $\IC_p$ is a separable complete valued field (so it  contains $\IQ_p$) which is algebraically closed (thus contains an algebraic closure of $\IQ_p$). 
  The valued field $\IC_p$ is  complete but not spherically  complete (see \cite[Section~3.4]{MR1760253}), hence not  locally compact.
  Notice that there are models of $\ZFA$  (see \cite{Lau62}) and of $\ZF$ (see \cite{Hodges76}) in which $\IQ$ has a (non well orderable hence non countable) algebraic closure $\IL$ without any non trivial absolute value on $\IL$.
 \end{remark} 

\subsection{Semi-normed vector spaces over a valued field}
Given a vector space $E$ over a valued field $(\IF,|.|)$, a {\em semi-norm} on $E$ is a mapping 
$N: E \to \IR_+$ satisfying for every $x,y \in E$ and  $\lambda \in \IF$ the properties  $N(\lambda.x)=|\lambda|N(x)$ and  $N(x+y) \le N(x)+N(y)$. 
For a ultrametric valued field $(\IF,|.|)$, the semi-norm $N$ is {\em ultrametric} if the semi-metric associated to $N$ is ultrametric, equivalently if for every $x,y \in E$, $N(x+y) \le \max(N(x), N(y))$. 
 
 \medskip

\section{In $\ZFA$, $\forall^{\Prime}p \; \MC(p)$ implies the whole Ingleton statement} \label{sec:mcp_to_ing}
\subsection{The Hahn-Banach axiom and the Ingleton axiom}
\subsubsection{$\AC$ implies the Hahn-Banach statement}
The following statement is a consequence of $\AC$: \\
$\HB$: (Hahn-Banach statement) {\em Given a $\IR$-vector space $E$, a semi-norm $N:E \to \IR_+$, a vector subspace  $V$ of $E$ and a linear form 
 $f: V \to \IR$ such that for every $x \in V$, $|f(x)|_{\IR} \le N(x)$, there exists a linear form  $\tilde f: E \to \IR$ extending 
 $f$ such that for every  $x \in E$, $|\tilde f(x)|_{\IR} \le N(x).$}

The usual proof of  $\HB$ can be obtained by transfinitely iterating the following Lemma  (for example using Zorn's lemma or a transfinite recursion on ordinals and the Axiom of Choice). 
\begin{lemma}[Hahn-Banach, 1932, ``one step''] \label{lem:one-stepHB}
Let $E$ be a $\IR$-vector space, let  $N: E \to \IR_+$ be a semi-norm on  $E$, let $V$ be a vector subspace of  $E$ and let  
$f: V \to \IR$ be a linear form such that $|f|_{\IR} \le N_{\restriction V}$. For every  $a \in E \backslash V$, there exists a linear form  $\tilde f: V + \IR.a \to \IR$ extending  $f$ such that  $|\tilde f|_{\IR} \le N_{\restriction V + \IR.a}$.
\end{lemma}

\begin{remark}
 In set-theory   without the axiom of choice $\ZFA$:
 \begin{enumerate}
   \item $\AC \Rightarrow \BPI \Rightarrow \RL \Rightarrow \HB$: 
see Jech's book \cite{Je73} or Howard and Rubin's book \cite{Ho-Ru}.
\item The implications  $\AC \Rightarrow \BPI $ and $\RL \Rightarrow \HB$ are not reversible in $\ZF$: see \cite{Je73} or \cite{Ho-Ru} for the first implication and \cite[Remark~9]{Mo-Rado} for the second one. 
It is known that $\RL$ does not imply $\BPI$ in $\ZFA$ (\cite{How84}) but the question ``Does $\RL$ imply 
$\BPI$?'' is open in
 $\ZF$.
\item $\HB \Rightarrow \text{``The Hausdorff-Banach-Tarski'' paradox}$ (see \cite{Pawli}). In set theory 
$\ZFA$, this implication is not reversible since the statement ``$\IR$ is well orderable'' implies the Hausdorff-Banach-Tarski paradox but does not imply $\HB$. 
 \end{enumerate}
\end{remark}

\subsubsection{$\AC$ implies Ingleton's statement}
The following ``one-step'' result is the ultrametric counterpart of Lemma~\ref{lem:one-stepHB}:
\begin{lemma}[Ingleton, 1952, ``one step'']  \label{lem:one-stepIng}
Let $E$ be a vector space over a spherically complete ultrametric valued field $(\IF,|.|)$, let $N: E \to \IR_+$ be a ultrametric semi-norm,  let $V$ be a vector subspace of  $E$ and let  $f: V \to \IF$ be a linear form such that $|f| \le N_{\restriction V}$. 
If  $a \in E \backslash V$, then there exists  a linear form  $\tilde f: V + \IF.a \to \IF$ extending  $f$ such that $|\tilde f| \le N_{\restriction V + \IF.a}$.
\end{lemma}

\begin{remark} Notice that both one-step Lemmas~\ref{lem:one-stepHB} and \ref{lem:one-stepIng} extend to a ``finite number of steps''.
\end{remark}

  For each spherically complete ultrametric valued field $(\IF,|.|)$, the Axiom of Choice implies the following statement (see~\cite{Ing52}):\\
 $\I_{(\IF,|.|)}$ (Ingleton's statement): {\em ``Let $E$ be a $\IF$-vector space, let $N: E \to \IR_+$ be a ultrametric semi-norm,  let $V$ be a vector subspace of  $E$ and let  $f: V \to \IF$ be a linear form such that $|f| \le N_{\restriction V}$. 
Then there exists a linear form    $\tilde f: E \to \IF$ extending  $f$ such that $|\tilde f| \le N$.''}

\medskip
 We shall  show that in set theory $\ZFA$, the  ``full Ingleton theorem'' + $\HB$ does not imply $\AC$ (unless $\ZFA$ is inconsistent). This answers in $\ZFA$ a question raised by van Rooij (see~\cite{Rooij}).

\subsection{A model of  $\ZFA \mathbf + \neg \AC$ with ``multiple choices''}
 Levy (1962) built a model $\mathcal N6$ of $\ZFA$ (see \cite[p.~185]{Ho-Ru}, \cite{Levy62}) in which there exists a sequence   $(F_n)_{n \in \IN}$ of finite sets such that for every $n \in \IN$, $\#F_n=n+1$ and  $\prod_{n \in \IN}F_n= \varnothing$: such a model  does not satisfy  $\AC$.
However, Levy showed that this model satisfies the following {\em ``Multiple Choice''} axiom: \\
\noindent   $\MC$: (``Multiple Choice'') {\em ``For every family  $(A_i)_{i \in I}$ of non-empty sets, there exists a family  $(B_i)_{i \in I}$ of non-empty finite sets such that for every $i \in I$,  $B_i \subseteq A_i$.''}\\
and also the following refined  statement for each  prime number  $p \ge 2$:\\
 $\MC(p)$: {\em ``For every family $(A_i)_{i \in I}$ of nonempty sets, there exists a family $(B_i)_{i \in I}$ of finite  sets such that for every  $i \in I$, $B_i \subseteq A_i$ and  $\#B_i$ is not a multiple of  $p$.''}

 \begin{remark}
 In set-theory $\ZFA$, $\AC \Rightarrow \forall^{\Prime}p \; \MC(p) \Rightarrow \MC$ and none of these implications is reversible.  
 In set-theory $\ZF$ (without atoms), $\MC$ implies $\AC$. 
 \end{remark}

\subsection{$\forall^{\Prime}p \; \MC(p)$  implies the  ``Full Ingleton'' statement}
Given a valued field $(\IF,|.|)$, a $\IF$-vector space $E$ endowed with a $\IF$-semi-norm $N$, 
a {\em finite linear extender} on $E$ is  a mapping associating to each ordered pair  $(V,f)$ where $V$ is a proper vector subspace of $E$  and $f: V \to \IF$ is a linear form such that   $|f| \le N_{\restriction V}$, an ordered pair  $(V',f')$ such that  $V'$ is a vector subspace of  $E$ strictly including  $V$ such that the $\IF$-vector space $V'/V$ is finitely generated, and   a linear mapping $f': V' \to \IF$   extending $f$ with   $|f'| \le N_{\restriction V'}$.
We shall prove in  $\ZFA$ the following result:
\begin{theorem} \label{theo:mc-to-finite-extender}
Let  $(\IF,|.|)$ be a spherically complete ultrametric valued field or the usual valued field $\IR$. Let  $E$ be a  $\IF$-vector space endowed with a semi-norm $N$  which is assumed to be ultrametric if  $\IF \neq \IR$.
\begin{enumerate}
\item  \label{it:MC} If $\IF$ is the archimedean field $\IR$ or if   $\IF$ has characteristic zero and if  $|.|_{\restriction \IQ}$ is the trivial absolute value, then $\MC$ implies
a finite linear extender on $E$.
\item   \label{it:MCp} In the other cases, there exists a  prime number $p$ such that the field $\IF$ has characteristic $p$ or such that  the ultrametric valued field   $(\IF,|.|)$  has characteristic zero and 
 $|.|_{\restriction \IQ}$ is equivalent to $|.|_p$; in these cases,   $\MC(p)$ implies a finite linear extender on $E$. 
\end {enumerate}
 \end{theorem}
 \Pr \eqref{it:MC}  With  $\MC$, let   $\Phi$ be a mapping associating to each non-empty subset $X$ of  $E \cup \IF^E$ a finite non-empty subset of  $X$.
Given a proper vector subspace $V$ of  $E$ and a linear form   $f: V \to \IF$ satisfying  $|f| \le N_{\restriction V}$,  
let $F:= \Phi(E \backslash V)$ and let $V_F:=\spanv(V \cup F)$. Using  Hahn-Banach's one-step 
Lemma~\ref{lem:one-stepHB} (for $\IF=\IR$) or Ingleton's one-step Lemma~\ref{lem:one-stepIng} (otherwise),
 the set   $\mathcal G$ of linear forms   $g: V_F \to \IF$ extending $f$ such that  
 $|g| \le N_{\restriction V_F}$ is non-empty. 
 Let  $G:=\Phi(\mathcal G)$.  Consider the linear form  
 $\tilde f:= \frac{1}{\#G}\sum_{g \in G}g$ on $V_F$: then $\tilde f$ extends  $f$.  
If $\IF=\IR$, then  $|\#G|_{\IR}=\#G$ thus  
 $|\tilde f(x)|_{\IR}  = \frac{1}{\#G }|\sum_{g \in G}g(x)|_{\IR} \le \frac{1}{\#G }\sum_{g \in G}|g(x)|_{\IR} \le N_{\restriction V_F}(x)$. If  $\IF$ has  characteristic zero and the restriction $|.|_{\restriction \IQ}$ is  the trivial absolute value, then $|\#G|=1$ thus for every $x \in V_F$, 
 $|\tilde f(x)|= \frac{1}{|\#G|}|\sum_{g \in G}g(x)|=
 |\sum_{g \in G}g(x)|
 \le \max_{g \in G}|g(x)| \le N(x)$ whence $|\tilde f| \le N_{\restriction V_F}$.\\
 \eqref{it:MCp}  If the characteristic of the field $\IF$ is zero, then $\IF$ extends the field $\IQ$ of rational numbers. Since   $|.|_{\restriction \IQ}$ is non-trivial,  Ostrowski's theorem implies that the absolute value induced by $|.|$ on $\IQ$ is equivalent to the $p$-adic absolute value for some prime number $p$.   With  $\MC(p)$, let  $\Phi_p$ be a mapping associating to each non-empty subset $X$  of $\IF^E$ a finite subset $G$ of $X$ such that  $p$ does not divide   $\#G$. 
Let $G:= \Phi_p(\mathcal G)$: then $G$ is a finite subset of $\mathcal G$ such that $p$ does not divide $\#G$. 
Let $n:=\#G$. Then $|n|=1$:  in the first subcase, $n \in \IF_p \backslash \{0\} \subseteq \IF$ 
thus $|n|=1$; in the second subcase, $|n|_p=1$   because $p$ does not divide $n$, thus $|n|=1$. \\
 Now we consider the linear form  
 $\tilde f:= \frac{1}{n}\sum_{g \in G}g$: 
 this linear form extends $f$, and for every  $x \in  V_F$, 
 $|\tilde f(x)| = \frac{1}{|n|}|\sum_{g \in G}g(x)| =
 |\sum_{g \in G}g(x)| \le \max_{g \in G}(|g(x)|) \le N(x)$, whence $|\tilde f| \le N_{\restriction V_F}$.
\Rp

\begin{corollary} \label{coro:MCptoIng} 
\begin{enumerate}
  \item \label{it:MCtoIng0}  $\MC$ implies $\HB$ and $\I_{(\IF,|.|)}$ for every spherically complete ultrametric  field 
$(\IF,|.|)$ with zero characteristic such that $|.|_{\restriction \IQ}$ is trivial. 
\item \label{it:MCptoIngp} For every prime number $p$, $\MC(p)$ implies  $\I_{(\IF,|.|)}$ for every spherically complete ultrametric  field with characteristic $p$ or with zero characteristic  such that  
 $|.|_{\restriction \IQ}$ is equivalent with the absolute value $|.|_p$ of $\IQ$. 
\end{enumerate}
\end{corollary}
\Pr Proof by transfinite recursion using  Theorem~\ref{theo:mc-to-finite-extender}.
\Rp

\begin{remark} The implication $\MC \Rightarrow \HB$ was already proved in \cite{Mo-Rado} ($\MC$ implies Rado's Lemma which implies $\HB$). The implication $\MC \Rightarrow  \I_{(\IF,|.|)}$ was proved 
for every spherically complete ultrametric  field $(\IF,|.|)$ with null characteristic  such that 
 $|.|_{\restriction \IQ}$ is trivial (see \cite[Proposition~4.4]{Mo-LinExt}). 
\end{remark}

It follows that the ``full Ingleton theorem'' follows from $\forall^{\Prime}p \; \MC(p)$  in set theory $\ZFA$,
thus in $\ZFA$, the full Ingleton statement does not imply $\AC$.

\subsection{Some questions}
\begin{enumerate}
  \item Are there links in set-theory without choice between the statements $\I_{\IF}$ obtained for various spherically complete ultrametric valued fields $\IF$? 
  \item Given a prime number $p$:
  \subitem  -Denoting by $\IF_p$ the finite field with $p$ elements, 
  does $\I_{\IF_p}$ implies $\I_{\IF}$ for every spherically complete ultrametric valued field with characteristic $p$?
  \subitem -Does $\I_{\IQ_p}$ imply $\I_{(\IF,|.|)}$ for every spherically complete ultrametric valued field $(\IF,|.|)$ with null characteristic such that $|.|_{\restriction \IQ}$ is equivalent to $|.|_p$?
  \item Does the conjunction of the statements $\I_{\IQ_p}$ for $p$ prime number imply $\I_{(\IQ,|.|_{triv})}$ or $\HB$? 
  \item Given two different prime numbers $p$ and $q$, are the statements $\I_{\IQ_p}$ and $\I_{\IQ_q}$ equivalent?
  \end{enumerate}

\begin{remark}
 For each ultrametric spherically complete valued field $(\IF,|.|)$, the statement $\I_{(\IF,|.|)}$ 
is equivalent to the following one (see \cite{Mo-LinExt}):\\
{\em ``For every vector subspace $F$ of an ultrametric semi-normed $\IF$-vector
space $(E,N)$, there exists an isometric linear extender $T : BL(F, \IF) \to  BL(E, \IF)$.''}
Here, given a vector subspace $V$ of $E$,  $BL(V,\IF)$ denotes the set of linear bounded mappings from
 $V$ to $\IF$.
 \end{remark}

\section{In $\ZFA$, $\MC$ implies the Krein-Milman statement $\KM$} \label{sec:MC2KM}
\subsection{Facets of subsets of a real vector space}
Given a {\em real} vector space $E$, and two elements $a,b \in E$, we denote by $[a,b]$ the {\em segment} subset $\{t.a + (1-t).b : t \in [0,1]\}$ of $E$; if $a=b$, then $[a,b]$ is a singleton, else 
the segment $[a,b]$ is infinite (it is equipotent with $[0,1]$)  and we denote by $]a,b[$ the 
{\em strict segment}  $]a,b[:=\{t.a + (1-t).b : t \in ]0,1[\}$.  
A subset $C$ of $E$ is {\em convex} if for every $a,b \in C$, $[a,b] \subseteq C$. 
In particular, every segment is convex. Say that an element  $e$ of a convex subset $C$ of $E$ is
an  {\em extreme point} of $C$ is $C \backslash \{e\}$ is convex: this means that for every
distinct elements $a,b$ of $C$, if $e \in [a,b]$ then $e=a=b$.
Say that a subset $F$ of a  subset $X$ of $E$ is a {\em facet} of $X$ if $F$ is nonempty, and if for every distinct elements $a,b$ in $X$, if $[a,b] \subseteq X$ and 
$]a,b[$ meets $F$ then $[a,b] \subseteq F$.

\begin{remark} A facet of a convex subset $X$ may be not convex: for example $\{0,1\}$ is a non-convex facet of the subset $[0,1]$ of the one-dimensional vector space $\IR$.
\end{remark}

\begin{remark} Given an element $e$ of a  convex subset $C$ of $E$, $e$ is an  extreme point  of  $C$  if and only if  $\{e\}$ is a facet of $C$.
\end{remark}

\begin{remark} \label{rem:facets}
Given a subset $X$ of a real vector space $E$,
\begin{enumerate}
\item \label{it:facet1} Every facet of a facet of $X$ is a facet of $X$;
\item \label{it:facet2}   A nonempty set which is an  intersection of facets of $X$ is a facet of $X$;
\item \label{it:facet3}  A union of a nonempty set of facets of $X$ is a facet of $X$;
\item \label{it:facet4}  if $X$ is nonempty  then $X$ is a facet of $X$.
\end{enumerate}
\end{remark}

Given a  subset $X$ of a real vector space $E$, 
a mapping $f: X \to \IR$ is said to be {\em convex} if for every $a, b \in X$ and every
$s \in [0,1]$, if $s.a+(1-s).b \in X$ then $f(s.a+(1-s).b) \le s f(a)+(1-s) f(b)$. The
mapping $f: X \to \IR$ is said to be {\em concave} if $-f$ is convex.  

\begin{lemma} \label{lem:facet-max}
Given a real vector space $E$, a subset $X$ of $E$ and a convex mapping 
$f: X \to \IR$ which is upper bounded and attains its least upper bound $M$, then the subset 
$F:=\{x \in X : f(x)=M\}$ is a facet of $X$.
\end{lemma}
\Pr Assume that there are two distinct elements $a,b \in X$ such that  $]a,b[$ meets $F$ and 
$[a,b] \subseteq X$. 
Let us show that $[a,b] \subseteq F$. Since $]a,b[$ meets $F$, let  $u \in ]a,b[$ 
such that $f(u) =M$: then there exists some real number $s$ such that 
$0<s<1$ and $u=s.a+(1-s).b$. By convexity of $f$, 
$M=f(s.a+(1-s).b) \le sf(a)+(1-s)f(b)$; if $f(a)<M$ or $f(b)<M$, 
then $sf(a)+(1-s)f(b) <M$ which is contradictory. Thus $f(a)=f(b)=M$. Now if $x \in ]a,b[$ and $x \neq u$,
then $u \in ]a,x[$ or $u \in ]x,b[$ whence  $f(x)=M$ by convexity of $f$. Thus $[a,b] \subseteq F$.
\Rp

Given a set $X$, a set $\mathcal C$ of subsets of $X$ is said to satisfy the {\em finite intersection property} if  every finite nonempty subset of $\mathcal C$ has a nonempty intersection. 
A subset $X$ of a real  linear topological vector space $E$ is {\em convex-compact} 
(see  \cite[p.~135]{Lu}) if for every family $\mathcal C$ of closed {\em convex} subsets of $E$ such that $\{C \cap X : C \in \mathcal C\}$
satisfies the finite intersection property,  $\cap \mathcal C \cap X$ is nonempty.

\begin{lemma} \label{lem:facet}
Given a real  topological vector space $E$, a convex-compact subset $X$ of $E$ and a continuous concave mapping $f: E \to \IR$ which is not constant on $X$, then $f$ is upper  bounded on $X$ and
attains its supremum  $M$ on $X$; moreover, if $f$ is linear, then the  closed subset  
$\{x \in X : f(x)=M\}$ of $X$ is a facet of $X$. 
\end{lemma}
\Pr  For every real number $\lambda \in f[X]$, let 
$C_{\lambda}:= \{x \in E : f(x) \ge \lambda\}$; since $f$ is continuous and concave,
$C_{\lambda}$ is a closed convex subset of $E$; moreover, 
since $\lambda \in f[X]$, $X \cap C_{\lambda} \neq \varnothing$. 
Let $\mathcal C := \{C_{\lambda} : \lambda \in f[X]\}$. Then 
$\{C \cap X : C \in \mathcal C\}$ satisfies the finite intersection property, 
thus $F:=\cap \mathcal C \cap X =\{x \in X : f(x)=M\}$ is nonempty; let $x_0 \in F$: then 
$M:=f(x_0)=\sup_X f$. Moreover, if $f$ is linear, then $f$  is convex, thus Lemma~\ref{lem:facet-max} implies that the  closed  subset 
$F$ of $X$  is a facet of $X$.
\Rp

\begin{remark}
The statement $\HB$ is equivalent to the following statement: 
\par {\em Given a real Hausdorff locally convex  topological vector space $E$, for every non null element   $x \in E$, there exists a continuous linear functional $f:E \to \IR$ such that 
$f(x) =1$.}
\end{remark}
\Pr $\Rightarrow$ Given some real Hausdorff locally convex  topological vector space $E$ and some non-null element $x$ of $E$, consider some convex open neighbourhood $V$ of $0_E$ such that $x \notin V$. Then, Theorem~2 of \cite{Do-Mo} implies the existence of a linear 
functional $f: E \to \IR$ such that $\sup_V f < f(x)$. Since $f: E \to \IR$ is linear and bounded on a neighbourhood of $0_E$, $f$ is continuous.  
For the converse statement, see \cite[Lemma~5]{Fo-Mo}.
\Rp

\begin{remark}
It follows that given a  subset $K$ of a real locally convex Hausdorff topological vector space $E$, if $K$ is convex-compact in $E$, then 
$\HB$ implies that closed  minimal facets of $K$ are singletons  of $K$.
\end{remark}
\Pr Given a closed  facet $F$ of $K$ which is not  a singleton, consider with $\HB$ some continuous linear functional $f: E \to \IR$ which is not constant on $F$; then, the closed subset $F$ of $K$ is convex-compact, thus 
 $f$ attains its upper bound $M$ on $F$; using Lemma~\ref{lem:facet}, the set
$\{x \in F : f(x)=M\}$ is a proper facet of the  set $F$, thus the facet $F$ of $K$ is not minimal. 
\Rp

\subsection{Trees  of subsets of a set}
\subsubsection{Trees}
A {\em tree} is a partially ordered set $(T,\preceq)$  with a smallest element $r$ (called the {\em root of the tree}) such that for every $x \in T$, the interval $\{y \in T : y \preceq x\}$ is well ordered.
Thus a tree is a well founded partially ordered set $(T,\preceq)$ such that for every $x \in T$, the interval $\{y \in T : y \preceq x\}$ is linearly ordered. 

\begin{remark}  Given an element $x$ of a tree $(T,\preceq)$, the set $\{y \in T : y \preceq x\}$ is the smallest  chain of $T$ containing $x$ and the root of $T$. 
\end{remark}

Given an element $x$ of a tree $(T,\preceq)$, we denote by $x^+$ the set of elements $y \in T$ which {\em cover} $x$ {\em i.e.} such that $x \prec y$ (where $\prec$ is the strict-order associated to $\preceq$) and the interval $]x,y[$ is empty; elements of $x^+$ are called {\em successors} of $x$. Given elements $x,y \in T$, if $y$ is a successor of $x$, then $x$ is said to be a {\em predecessor} of $y$. Every element of $X$ has at most one predecessor in the tree $T$.
A {\em leaf} of the tree $(T,\preceq)$ is a maximal element of the {\em poset} 
$(T,\preceq)$. Thus an element $x$ of a tree $T$ is a leaf of $T$ if $x$ has no successor in $T$. 
 A {\em branch} of the  tree 
$(T,\preceq)$ is a maximal chain of  the {\em poset} $(T,\preceq)$: every branch has a  first element which is the root of the tree. 
Leaves of a tree correspond to greatest elements of branches of this tree.  A {\em subtree} of a tree $(T,\preceq)$ is a nonempty initial section of the {\em poset} $(T,\preceq)$ {\em i.e.} a subset of $T$ containing the root of $T$ and such that for each $x \in S$,
$\{y \in T : y \preceq x\} \subseteq S$.  Of course, each subtree of a tree $(T,\preceq)$ is a tree for the order induced by $\preceq$ on the subtree.

\begin{remark} In $\ZFA$, the Axiom of Choice implies that  every tree has a branch. The converse implication is also true since the existence of a branch in every tree implies the statement $\forall \kappa \DC_{\kappa}$ (form {1F} in \cite[p.~12]{Ho-Ru}) which in turn implies $\AC$ 
(\cite[Th.~8.2. p.~121]{Je73}). 
\end{remark}

Given a tree $(T,\preceq)$ with root $r$, since the {\em poset} $(T,\preceq)$ is well-founded, we may consider the ``rank function'' $\rk$ associating to each 
$x \in T$ the  ordinal such that $\rk(r)=0$ and, for each $x \in T \backslash \{r\}$, $\rk(x) = \sup \{\rk(y) +1 : y \prec  x\}$. 
 An element $y \in T$ has a predecessor if and only if the ordinal $\rk(y)$ is a successor ordinal.  
  Denoting by  $\On$ the collection of ordinals, we  consider the family $(L_{\alpha})_{\alpha \in \On}$ of {\em level sets} of the tree $T$, where  for each ordinal $\alpha$, $L_{\alpha}$ is the set of elements of $T$ with rank $\alpha$.  
  Thus $L_0=\{r\}$ where $r$ is the root of the tree $T$.  
  For each ordinal $\alpha$, we denote by $T_{\alpha}$ the set $\{x \in T : \rk(x)< \alpha\}$.  Thus $T_0=\varnothing$, and for each ordinal 
  $\alpha \ge 1$, 
$T_{\alpha}$ is a subtree of $T$. Moreover, for each ordinal $\alpha$, $T_{\alpha+1}=T_{\alpha} \cup L_{\alpha}$, and for each limit ordinal 
$\alpha >0$, $T_{\alpha}=\cup_{\beta \in\alpha}T_{\beta}$.

Notice that for every ordinal $\alpha$,  if $L_{\alpha}=\varnothing$, then for every $\beta \ge \alpha$, $L_{\beta}=\varnothing$ and $T_{\beta}=T_{\alpha}$. Since $T$ is a set, the axiom schema of replacement  implies that  there exists an ordinal $\alpha$ such that 
$T_{\alpha+1}=T_{\alpha}$ (whence $L_{\alpha}=\varnothing$). The first ordinal $\alpha$ such that $L_{\alpha}$ is empty is called the {\em rank of the tree $T$}.

\subsubsection{Trees   of subsets of  a set $X$}
Given a set $X$, a  {\em tree of subsets}  of $X$ is a subset $T$ of the {\em poset} $(\mathcal P(X) \backslash \{\varnothing\}),\supseteq)$, with smallest  element $X$,
such that the induced {\em poset} on $T$ is a   tree, and such that for every $x \in T$, successors of $x$ are {\em pairwise disjoint} subsets of $X$. 
It follows that for every $x, y \in T$, if $x$ and $y$ are not comparable for the inclusion relation then $x \cap y=\varnothing$. 

\begin{remark} Given a tree $T$ of subsets of a set $X$ with rank $\alpha$, and denoting for each $\beta \in \alpha$ by $L_{\beta}$ the level set $\{x \in T : \rk(x)=\beta\}$, and by $F_{\beta}$ the subset $\cup L_{\beta}$ of $X$, then the family $(F_{\beta})_{\beta \in \alpha}$ is descending 
{\em i.e.} for every $s,t \in \alpha$ such that $s \in t$, $F_t \subseteq F_s$.  
\end{remark}
\Pr Given $b \in F_t$, there exists $y \in L_t$ such that $b \in y$. Let $x$ be the element of $L_s$ such that $x \preceq y$; then $y \subseteq x$
thus $b \in x \subseteq F_s$.
\Rp

\begin{lemma} \label{lem:branch-tree} If $T$ is a tree of nonempty closed subsets of a compact topological space $X$, and if each level of $T$ is finite, then  $T$ has at least one branch. 
\end{lemma}
\Pr Let $\alpha$ be the rank of the tree $T$.  
For each $\beta < \alpha$, let $F_{\beta}$ be the (nonempty) closed subset $\cup L_{\beta}$. Then the family $(F_{\beta})_{\beta \in \alpha}$ is decreasing thus the set $\cap_{\beta < \alpha}F_{\beta}$ is nonempty: let $x \in \cap_{\beta < \alpha}F_{\beta}$;  let $b:=\{F \in T : x \in F\}$;
since  elements of $T$ are pairwise disjoint,  $b$ is a branch of $T$. 
\Rp

\subsubsection{Extending branches of a tree of subsets}
Given a set $X$ and a tree $T$ of subsets of $X$, 
given a branch $b$ of $T$ with no leaf, if $\cap b \neq \varnothing$,  
and if $F$ is a set of pairwise disjoint nonempty subsets of $\cap b$, 
then $b \cup F$  is a tree of subsets of $X$: we call it the tree 
{\em obtained by extending the branch $b$ with successors in  $F$}.
 Given a nonempty set $B$ of branches of $T$, if for each $b \in B$, $F_b$ is a nonempty set of nonempty pairwise disjoint subsets of $\cap b$, 
then  $\cup B \cup \cup_{b \in B} F_b$ is a tree of subsets of $X$ which is called the
 {\em tree obtained from the subtree $\cup B$ by extending each $b \in B$ by the set of successors $F_b$}.  

\subsubsection{Dynasties of  trees of subsets}
A {\em dynasty  of trees} of subsets of $X$ is a collection $(T_{\alpha})_{\alpha \in \On}$ of trees of subsets of $X$, such  for each 
$\alpha \in \On$, there exists a set $B$ of branches of the tree   $S:=\cup_{\beta< \alpha}T_{\beta}$ such that either $T_{\alpha}$ is the subtree $\cup B$ 
of $S$, or $T_{\alpha}$   is obtained  from the tree $\cup B$  by extending each   branch  $b \in B$ by a nonempty set of successors.

\subsection{The axiom of ``Multiple Choices'' $\MC$ implies the Krein-Milman axiom}

\begin{lemma} \label{lem:pairwise-disjoint}  Given a real  vector space $E$ and  a nonempty subset $X$ of $E$, there is a mapping associating to each nonempty finite set $F$ of convex facets of $X$  a nonempty finite set $G$ of {\em pairwise disjoint} convex facets of $X$ such that for each $y \in G$, there exists $x \in F$ such that $y \subseteq x$.
\end{lemma}
\Pr Let $\mathcal F$ the set of convex facets of $X$. Given a nonempty finite subset $F$ of $\mathcal F$, let $F_1$ be the set of maximal subsets $J$ of $F$ such that $\cap J$ is nonempty; then $G:=\{ \cap J : J \in  F_1\}$ is a set of convex facets of $X$ such that for each $C \in G$, there exists $C' \in F$ satisfying $C \subseteq C'$; elements of $G$ are pairwise disjoint because elements of $F_1$ are maximal subsets of $F$ with nonempty intersection. 
\Rp

\begin{theorem} \label{theo:mc2hb}
In $\ZFA$, $\MC$ implies $\KM$.
\end{theorem}
\Pr Let $E$ be a real locally convex Hausdorff topological vector space and let $K$ be a nonempty convex compact subset of $E$. Using $\MC$, 
let $\Phi$ be a mapping associating to each nonempty set $A$ of subsets of $K$ a nonempty finite subset of $A$.
We shall build by recursion a dynasty   $(T_{\alpha})_{\alpha \in \On}$  of trees of closed convex facets of $K$, such that  for each $\alpha \in \On$, the tree $T_{\alpha}$ has finite levels. 
 The tree $T_0$ is the singleton $\{K\}$ (here $K$ is a  facet of $K$).  Given an ordinal $\alpha >0$ such that the family $(T_{\beta})_{\beta \in \alpha}$ has been defined, and denoting by $S$ the tree $\cup_{\beta < \alpha}T_{\beta}$, Lemma~\ref{lem:branch-tree} implies that the set of branches of the tree $S$ is nonempty (because $S$ is a tree of closed subsets of a compact space and each level of $S$ is finite). Using $\Phi$,  we consider a nonempty finite set $B$ of branches of $S$, and we   distinguish two  cases according to whether there exists $b \in B$  such that
  $\cap b$  is a minimal facet of $K$  or not. 
 In the first case and if  there exists $b \in B$ such that $b$ has a leaf, then $T_{\alpha}:=\cup \{b \in B : b \text{ has a leaf in } S\}$.
 If  no branch $b \in B$ has a leaf, then, for each $b \in B$, $F_b:=\cap b$ is a minimal facet of $K$  and we define 
 $T_{\alpha}= B[({F_b})_{b \in B}]$. 
 In the second case, we consider, with $\Phi$ and Lemma~\ref{lem:pairwise-disjoint}, for each $b \in B$, a nonempty finite set $F_b$ of {\em pairwise disjoint} convex proper facets  
 of $\cap b$, and we define the tree  $T_{\alpha}=B[({F_b})_{b \in B}]$ obtained from the subtree $\cup B$ by extending each $b \in B$ by the set $F_b$ of successors.
 Using the axiom schema of replacement, there exists an ordinal  $\alpha$  such that 
$T_{\alpha}=T_{\alpha+1}$, thus  $T_{\alpha}$ has a leaf which is a closed minimal convex facet of $K$: this facet corresponds to an extreme point of $K$. 
\Rp

\begin{corollary}
In $\ZFA$, $\MC$ implies that every compact  convex subset of a real Hausdorff locally convex topological vector space is the closed convex hull of the set of its extreme points.
\end{corollary}
\Pr Given  a nonempty compact  convex subset $K$ of a real Hausdorff locally convex topological vector space $E$, $\MC$ implies that the set $X$ of extreme points of $K$ is nonempty  (see Theorem~\ref{theo:mc2hb}); then using $\HB$ (which is a consequence of  $\MC$), it follows that $K$ is the closed convex hull of $X$. 
\Rp

\begin{remark} Consider the following  statement $\VKM$: ``Every nonempty convex-compact convex subset of a real locally convex topological vector space has an extreme point.'' (form 286 in \cite{Ho-Ru}).  It is known that $\HB + \VKM$ implies $\AC$ (\cite{B-F}), thus $\forall^{\Prime}p \; \MC_p$ does not imply $\VKM$. Also notice that $\VKM$  implies $\forall \kappa \AC^{\kappa}$ (see \cite{Fo-Mo}). However, it is an open question whether $\VKM$ implies $\AC$.  
\end{remark}

\bibliographystyle{abbrv}
\bibliography{Ingleton_AC.bbl}
\end{document}